\documentclass[11pt]{amsart}
\usepackage[T1]{fontenc}
\usepackage{amsmath}
\usepackage{amsthm}
\usepackage{amssymb}
\usepackage{amsfonts}
\usepackage{outlines}
\usepackage{dsfont}
\usepackage{stackengine}
\usepackage{comment}
\usepackage{cancel}
\usepackage{graphicx}
\usepackage{enumerate}
\usepackage{mathrsfs}
\usepackage{tikz-cd}
\usepackage{hyperref}
\usepackage{enumitem}
\usepackage{enumerate}
\usepackage{lipsum}
\usepackage[square,numbers]{natbib}
\usepackage[left = 3.5cm, right = 3.5cm, top = 4 cm, bottom = 4.5 cm]{geometry}

\newtheorem{thm}{Theorem} [section]

\newtheorem{lemma}[thm]{Lemma}

\newtheorem{remark}[thm]{Remark}
\newtheorem{prop}[thm]{Proposition}

\newcommand{\mres}{\mathbin{\vrule height 1.5ex depth 0pt width 0.13ex\vrule height 0.13ex depth 0pt width 1.3ex}}

\newcommand*{\RangeX}{%
 \mathsf{X}\mkern-7.5mu\mathsf{X}%
}

\newenvironment{claim}[1]{\par\noindent\textbf{Claim:}\space#1}{}
\newenvironment{claimproof}[1]{\par\noindent\underline{Proof:}\space#1}{\hfill $\square$}

\title[LMCF of surfaces with mean curvature bound]{Lagrangian mean curvature flow of surfaces with mean curvature bound}
\author{Sourav Ghosh}
\date{Sep 30, 2024}

\address{
Department of Mathematics, University of Notre Dame, Notre Dame, IN 46556}
\email{sourav.ghosh.1@nd.edu}

\begin{document}

\begin{abstract}
Let $L_t$ be a zero Maslov Lagrangian mean curvature flow in $\mathbb{C}^2.$ We show that if the mean curvature stays uniformly bounded along the flow, then the tangent flow at a singular point is unique i.e. the limit of the parabolic rescalings does not depend on the chosen sequence of rescalings.
\end{abstract}

\maketitle 

\tableofcontents

\section{Introduction}\label{}

Let $M$ be a smooth surface. A family of embeddings $F: [0,T) \times M \to \mathbb{R}^n$ satisfies the mean curvature flow if they satisfy the equation
$$\Big(\frac{\partial F}{\partial t}\Big)^{\perp} = H.$$
Here $H = H_{M_t}$ denotes the mean curvature vector of the embedded submanifolds $M_t := F_t(M) \subseteq \mathbb{R}^n$ and $\perp$ denotes the projection to the normal bundle. Submanifolds evolving by mean curvature flow often develop singularities in finite time $T < \infty$.

A new side of mean curvature flow was uncovered by Smoczyk \cite{smoczyk1996canonical} when he proved that the mean curvature flow preserves Lagrangian submanifolds of Calabi-Yau manifolds, which is known as Lagrangian mean curvature flow. Lagrangian submanifolds of Kahler manifolds are those for which the symplectic form vanishes. Since mean curvature flow is the gradient descent of volume, a natural hope is that the Lagrangian mean curvature flow enjoys long-time existence and convergence to a minimal Lagrangian representative – a special Lagrangian. 

Just as in mean curvature flow in general, the Lagrangian mean curvature flow forms singularities. In fact, by a result of Neves \cite{neves2013finite} in the langrangian mean curvature flow, singularities are unavoidable along even if the initial Lagrangian is a small Hamiltonian perturbation of a special Lagrangian. Along the mean curvature flow of zero-Maslov Lagrangians, all tangent flows at singularities are given by unions of minimal Lagrangian cones, according to Neves \cite{neves2007singularities}. In particular all such tangent flows are either singular, or have higher multiplicity. 

The simplest kind of singularities, called neck pinches in the two dimensional case has been studied well in Lotay-Schulze-Sz\'ekelyhidi \cite{lotay2022neck}. They have shown that if one tangent flow is given by a union of two multiplicity one transverse planes, then it is unique. Using that they have shown that the flow can be continued past the singular time if the singularities are modelled on two transverse, multiplicity one planes. In \cite{lotay2024ancient} by the same authors, some progress were made in order to understand what we can say about the flow in the presence of such singularities given by two static planes meeting along a line.

Uniqueness of blowups is perhaps the most fundamental question that one can ask about singularities and implies regularity of the singular set. There are many results in mean curvature flow when the tangent flow is multiplicity one. It is worth mentioning some related uniqueness results. By a monotonicity formula of Huisken \cite{huisken1990asymptotic}, and an argument of Ilmanen \cite{ilmanen1995singularities} and White \cite{white2005local}, tangent flows are shrinkers, i.e., self-similar solutions of MCF that evolve by rescaling. There are various uniqueness results depending on the structure of the tangent flow. When a tangent flow of a mean curvature is compact, then Schulze \cite{schulze2014uniqueness} proved that it is the unique tangent flow. Colding-Minicozzi \cite{colding2015uniqueness} proved uniqueness of the tangent flow when a tangent flow is a generalized cylinder. Chodosh-Schulze \cite{chodosh2021uniqueness} proved uniqueness of the tangent flow when a tangent flow is smooth and asymptotically conical. The first uniqueness result for tangent flows given by non-smooth shrinkers came in \cite{lotay2022neck}, in the setting of Lagrangian mean curvature flow. 
Also, Stolarski \cite{stolarski2023structure} proved that for flows with bounded mean curvature, if a tangent flow is given by a regular multiplicity one stationary cone, then it is the unique tangent flow. But all the above mentioned results are applicable for multiplicity one. So far, we hardly know anything about tangent flows of higher multiplicity. 

The Multiplicity One Conjecture for mean curvature flows of surfaces in $\mathbb{R}^3$ says that any blow-up limit of such mean curvature flows has multiplicity one. Li-Wang \cite{li2019extension} proved this result in a very special case, where the mean curvature stays uniformly bounded and later in \cite{li2022ilmanen} the same authors proved the result under the assumption that the mean curvature is of type-I. Recently, the conjecture has been solved by Balmer-Kleiner \cite{bamler2023multiplicity} in full generality. But again this is a result in codimension one flow of surfaces but for Lagrangian MCF of surfaces higher multiplicity is not ruled out.

Our main theorem is the following,

\begin{thm} {\label{thm 1.1}}
    Let $(L_t)_{0 \leq t < T}$ be a zero Maslov Lagrangian mean curvature flow in $\mathbb{C}^2$ with uniformly bounded mean curvature that develops a finite time singularity at space-time $(x_0,T)$. Then the tangent flow at $(x_0,T)$ is unique.
\end{thm}

Theorem \ref{thm 1.1} is notable as little is known in general about tangent flows of higher multiplicity for mean curvature flow. The uniqueness result here applies not only to non-smooth minimal cones which may arise as tangent flows in the Lagrangian setting but also cones with higher multiplicity, for example, a plane with multiplicity two, although our proof only works under the assumption of a uniform mean curvature bound.

\subsection{Outline} Using basic techniques of Geometric measure theory we can show that if we assume a uniform mean curvature bound then for any sequence of time converging to the singular time, there exists a subsequence which converges to a varifold and a current at the final time slice. Furthermore, we can show that the varifold does not depend on the sequence chosen. 

We will pick one such limiting current and prove that the limiting current has an almost minimizing property using the Lagrangian condition. We will see that to analyze the tangent flow of our flow, it is enough to consider all the static flows of the tangent cones of the varifold associated to our almost minimizing limiting current. This is the main technical part of our result. A similar result has been proved by \cite{stolarski2023structure} for the general mean curvature flow for limiting varifolds. But for our purpose, we need a slightly stronger result from what has been proved in that paper. In this paper, we prove the result in the Lagrangian setting, using a different argument for the varifold associated with our limiting current.

Even though our limiting current may not be unique, we will show using the almost minimizing property that for any two different limiting currents the static flows of the tangent cones of the varifold associated to them are the same. We will mainly focus on smooth flows in $\mathbb{C}^2$, but this part of our argument will work for Lagrangian mean curvature flows of arbitrary dimension with mean curvature bounds. 

Then in the next step we show that there is a one to one correspondence between the tangent flow of the original flow and the tangent cone of any limiting current. So, to prove Theorem \ref{thm 1.1}, it will suffice to prove that the tangent cone of the limiting current is unique. As we are in the two-dimensional setting, White's epiperimetric type inequality [\citealp{white1983tangent}, \citealp{de2017uniqueness}] will help us to prove the last statement. This part of our result, the way we prove it, only works in $\mathbb{C}^2$.

\subsection{Acknowledgments}
I would like to thank my advisor, Professor G\'abor Sz\'ekelyhidi, for introducing me to this circle of ideas and his constant support and encouragement throughout the project. I like to thank Professor Nicholas Edelen for many helpful discussions at various stages. I also thank Dr Maxwell Stolarski for the helpful comments.
\section{Preliminaries}\label{}

\subsection{Preliminaries from Geometric measure theory}
Here we will briefly discuss the notion of varifold and various related concepts; further details can be found in Simon’s book \cite{simon1983lectures}. 

We will denote the Hausdorff measure of dimension $n$ by $\mathcal{H}^n$. For an open set $U \in \mathbb{R}^n$, let $\mathcal{G}_k(U)$ denote the Grassmanian bundle of $k$-dimensional hyperplanes over $U$. A $k$ varifold in $U$ is a positive Radon measure on $\mathcal{G}_k(U)$. The topology on the space of $k$-varifolds $\mathcal{V}_k(U)$ is given by the weak* topology i.e. a net $\{V_i\}_{i\in I} \subset \mathcal{V}_k(U)$ converges to $V$ if
$$\int_{\mathcal{G}_k(U)} f(x,S) dV_i (x,S) \to \int_{\mathcal{G}_k(U)} f(x,S) dV (x,S)$$
for all $f \in C_c(\mathcal{G}_k(U)).$

If $\phi : U \to U'$ is a diffeomorphism and $V \in \mathcal{V}_k(U)$, we define $\phi_{\#} V \in \mathcal{V}_k(U')$ by the following formula,
$$\phi_{\#} V (f) = \int_{\mathcal{G}_k(U)} f \big(\phi(x),d_x\phi(S)\big) J\phi(x,S) dV(x,S),$$
where,
$$J\phi(x,S) = \Big( det \big( (d_x \phi |_\omega)^t \circ (d_x \phi |_\omega) \big)  \Big)^{1/2},$$
and $f \in C_c(\mathcal{G}_k(U')).$
Given a compactly supported, smooth vector-field $X$ on $U$ let $\phi_t$ denote the flow of $X$; the first variation of $V$ is given by
$$\delta V (X) = \frac{d}{dt} {\Big|}_{0} \big(\phi_t\big)_{\#} \mu_V (U),$$
where $\Pi$ is the canonical projection from $\mathcal{G}_k(U)$ to $U$ and  $\mu_V(A) = V(\Pi^{-1} (A)).$

If $V$ has locally bounded first variation, then by Riesz Representation the total variation measure of $\delta V$ exists. We say that $V$ has generalized mean curvature $H \in L^1_{loc}(U)$ if there exists a Borel function
$H : U \to \mathbb{R}^n$ with $H \in L^1_{loc}(U, d \mu_V)$ such that the first variation $\delta V$ satisfies
$$\delta V(X) = - \int H.X \;d \mu_V.$$

Given a $k$ rectifiable set $S \subset U$ and a non-negative function $\theta \in L^1_{loc}(S, \mathcal{H}^k \mres S)$ we define the $k$ varifold $v(S,\theta)$ by
$$v(S,\theta)(f) = \int_S f(x,T_x S) \;\theta(x) \;d \mathcal{H}^k(x),$$
where $T_x S$ denotes the approximate tangent space of $S$ at $x$ which exists $\mathcal{H}^k \mres S$ almost everywhere. $V$ is called an integral $k$ varifold if $V = v(S, \theta)$ for some $S$ and $\theta$ with $\theta$ taking non-negative integer values $\mathcal{H}^k \mres S$ almost everywhere.

We have compactness theorem (Allard compactness theorem) in the set of integral varifolds. Suppose $V_i \subset \mathcal{V}_n(U),$ a sequence of integral varifold and
$$\sup\limits_{i} \big(V_i(W) + \|\delta V_i\| (W)) \big) < \infty  \hspace{0.5 cm} \forall \;W \Subset U,$$ 
then there exists a subsequence $i_k$ and an integer multiplicity varifold $V$ such that $V_{i_k} \rightharpoonup V$. 

By compactness theorem for varifolds, for any integer rectifiable varifold $V$ in $U$ with bounded generalized mean curvature, which enjoys a uniform bound for $\mu_V(K)$, whenever $K$ is compact, for any $x_0 \in U$ and any sequence $\lambda_i \to \infty,$ there exists a subsequence $\lambda_{i_j}$ and an integer rectifiable varifold $C$ such that  
$(F_{x_0, \lambda_{i_j}})_{\#} V$ converges to $C$ in the sense of varifold, where
$$F_{x_0, \lambda}: \mathbb{R}^n \to \mathbb{R}^n \hspace{1 cm} F_{x_0, \lambda}(x) := \lambda(x-x_0).$$
Any $C$ that arises as such a limit is called a tangent cone of $V$ at $x_0$.

Now let us discuss some basic current related concepts; again for further details one can look at \cite{simon1983lectures}. 

An $n$ dimensional current in an open set $U \in \mathbb{R}^n$ is a continuous linear functional on the set of smooth differential $n$ form which has compactly support. The topology on the set of $n$ dimensional currents $\mathcal{D}_n(U)$ is given by the weak* topology i.e. $T_i \rightharpoonup T$ iff 
$$T_i(\omega) = T(\omega) \hspace{0.5 cm}\forall \omega \in \mathcal{D}_n(U). $$ 
The mass of the current $T$ is defined by 
$$M(T) = \sup \limits_{\|\omega\| \leq 1, \;\omega \in \mathcal{D}_n(U)} T(\omega),$$
where, 
$$\|\omega\|^2 = \sup \limits_{x\in U}\langle\omega(x),\omega(x)\rangle.$$
More generally for any open $W \subset U$, 
$$M_W(T) = \sup \limits_{\|\omega\| \leq 1, \;\omega \in \mathcal{D}_n(U),\; spt \omega \subset W} T(\omega).$$
If $M_W(T) < \infty$ for every $W \Subset U,$ then by the Riesz Representation theorem, we can associate a radon measure associated to the current $T$, we will denote it by $\mu_T.$ Given $T \in \mathcal{D}_n(U)$ the support, spt $T$ of $T$ to be the relatively closed subset of $U$ is defined by
$$spt T = U \setminus \cup W,$$
where the union is over all open sets $W \Subset U$ such that $T(\omega) = 0$ whenever spt $\omega \in U$. For a smooth map $f$ such that $f|_{spt T}$ is proper, the pushforward of a current is defined by 
$$f_\#T(\omega) = T(\eta f^\# \omega),$$
where $\eta$ is a compactly supported smooth function, identically equal to one in a neighborhood of the compact set spt$T \;\cap$ spt$f^\# \omega.$ We define the cone over current $T$, denoted $0 \;\RangeX \;T$ is defined by 
$$0 \;\RangeX \;T = h_\# \big([0,1] \times T \big),$$
where $h$ is defined by $(t,x) \mapsto tx$.

If $T \in \mathcal{D}_n(U)$, we say that $T$ is an integer multiplicity rectifiable $n$-current if for any $\omega \in \mathcal{D}^n(U),$
$$T(\omega) = \int_M \langle \omega(x), \eta(x) \rangle \;\theta(x) d\mathcal{H}^n(x),$$
where $M$ is an $\mathcal{H}^n$-measurable countably $n$-rectifiable subset of $U$, $\theta$ is a locally $\mathcal{H}^n$-integrable positive integer-valued function, and $\eta(x)$ can be expressed in the form $\tau_1 \wedge \tau_2 \wedge ... \wedge \tau_n$, where $\tau_1 \wedge \tau_2 \wedge ... \wedge \tau_n$ form an orthonormal basis for the approximate tangent space $T_x M$ and $\eta$ is $\mathcal{H}^n$ measurable and $\langle \;,\;\rangle$ is the natural dual pairing. 

We have compactness theorem (Federer and Fleming compactness theorem) in the set of integral currents. Suppose $T_i \subset \mathcal{D}_n(U),$ suppose $T_i, \partial T_i$ are integer multiplicity for each $i,$ and
$$\sup\limits_{i} \big(M_W (T_i) + M_W(\partial T_i) \big) < \infty  \hspace{0.5 cm} \forall \; W \Subset U,$$ 
then there exists a subsequence $i_k$ and an integer multiplicity current $T$ such that $T_{i_k} \rightharpoonup T$. 

Using the compactness theorem for currents, if a current does not have boundary and enjoys uniformly mass bound on each compact sets, one can define tangent cone of that current in a similar way we defined for varifold. Note that tangent cone of a current will have a current (and so an associated varifold) structure, while tangent cone of a varifold has only varifold structure.
 
\subsection{Preliminaries from Lagrangian mean curvature flow}

In this section, we will briefly discuss the basic definitions in Lagrangian mean curvature flow that will be useful in the article. 

Let us denote the standard complex structure on $\mathbb{C}^2$ by $J$ and the standard symplectic form by $\omega.$ We will consider $\Omega$ for the holomorphic volume form $dz_1 \wedge dz_2$. A (real) $2$-dimensional submanifold $L$ in $\mathbb{C}^2$ is said to be \textit{Lagrangian} if $\omega|_L = 0$. One can check that 
$$\Omega|_L = e^{i\theta}vol_L,$$
where $vol_L$ denotes the induced volume form of $L$ and $\theta$ is some function (possibly multiple-valued) called the \textit{Lagrangian angle}. When $\theta$ is a single valued function, $L$ is called \textit{zero-Maslov class}. 

For lagrangian manifolds, one has the following relation between the Lagrangian angle and the mean curvature,
$$H = J \nabla \theta.$$
Let $L_0$ be a lagrangian in $\mathbb{C}^2$ with uniformly bounded area ratios i.e. for some constants $C_0 , R_0$, we have the area bounds
$$\mathcal{H}^2\big(L_0 \cap B_R(0)\big) \leq C_0 R^2$$
for all $R \geq R_0$. Then one can show [for instance \citealp{neves2007singularities}, Appendix $B$] that for all $t \geq t_0,$ and $x \in \mathbb{C}^2,$ there is a constant $C$ depending on $C_0, |x|$ and $R_0 t_0^{-1/2}$ such that 
$$\mathcal{H}^2\big(L_t \cap B_R(x)\big) \leq C R^2$$
for all $R>0.$ Now suppose $L_0$ has uniformly bounded area ratios and we have a solution $(L_t)_{\{0 \leq t < T \}}$ to mean curvature flow for which the second fundamental form of $L_t$ is
bounded by a time dependent constant. Then using the argument used in \cite{smoczyk1996canonical}, one can show that the Lagrangian condition is preserved by the flow. 

For $\lambda > 0,$ we have the parabolic rescaling,
\begin{align*}
    \mathcal{D}_\lambda: \mathbb{C}^n \times \mathbb{R} \to \mathbb{C}^n \times \mathbb{R} \\
    (x,t) \to (\lambda x, \lambda^2t).
\end{align*}
If we denote the space time track of mean curvature flow by $\mathcal{L}$, then for a (Lagrangian) mean curvature flow $\mathcal{L}$, one can check that $\mathcal{D}_\lambda \mathcal{L}$ is again a (Lagrangian) mean curvature flow. Assume now that the solution to mean curvature flow develops a singularity at the point $(x_0,T)$ in space-time. For a positive sequence $\lambda_i \to \infty,$ it is well known that $\mathcal{D}_{\lambda_i}\big(\mathcal{L}-(x_0,T)\big)$ has a subsequence which converges weakly i.e. as a Brakke flow [\citealp{ilmanen1994elliptic}, Lemma $7.1$]. Such a flow is called a tangent flow at $(x_0,T)$ and it depends on the sequence $\lambda_i$ taken.

Let $(L_t)_{0\leq t < T}$ be a solution to Lagrangian mean curvature flow. Let $x_0 \in \mathbb{C}^n.$  We will consider the rescaled flow $L_\tau$ at $(x_0,T)$ which is defined by
$$L_\tau(x) = \frac{1}{(T-t)^{1/2}} \big(L_t(x) - x_0 \big), \hspace{0.5cm} \tau(t) = -log(T-t).$$
It is easy to check that the rescaled flow evolves with normal speed $H + \frac{1}{2}x^\perp.$

\section{Limiting current of the flow}\label{}

First we will show that if we assume the mean curvature is uniformly bounded throughout the flow, then we can extract a weak limit at time $t = T.$ This is true for any mean curvature flow, we do not need the Lagrangian condition to prove this statement. Similar type of results have been proved for varifolds in \cite{stolarski2023structure}, and for both varifold and currents in \cite{neves2007singularities}. 

\begin{lemma} {\label{lemma 3.1}}
    Let $(L_t)_{t \in [0,T)}$ be a solution to the mean curvature flow with uniformly bounded mean curvature. Then there exists a unique integer rectifiable varifold $L_T$ such that as $t \to T, L_t$ converges to $L_T$ in the sense of varifold and for any sequence $t_i \to T,$ there exists a subsequence (still denoted $t_i$) such that $L_{t_i}$ converges in the sense of current. 
\end{lemma}
    
\begin{proof}
    Let $t_i$ be a sequence converges to $T.$ The first variation formula yields, for any vector field $X$ supported in $B_r(0),$
    $$\delta L_{t_i} (X) = - \int_{L_{t_i} \cap B_r(0)}  \langle H, X \rangle d \mathcal{H}^n + \int_{\partial L_{t_i} \cap B_r(0)}  \langle X, \nu \rangle d \mathcal{H}^{n-1}$$
    where $\nu$ denotes the exterior unit normal. Hence, whenever the sup norm
    of $\|X\|_\infty$ satisfies $\|X\|_\infty \leq 1$, we have
    $$|\delta L_{t_i} (X)| \leq Cr \Big (\int_{L_{t_i} \cap B_r(0)} |H|^2 \Big)^\frac{1}{2}$$
    So by Allard compactness theorem for varifolds and Federer and Fleming compactness theorem for integral rectifiable current, we can conclude that $L_{t_i}$ has a subsequence which converges weakly to integral varifold and integral current respectively. 
    
    Now we will show the limiting varifold does not depend on the sequence $t_i \to T$. Let there are two sequences $t_{i_1}$ and $t_{i_2}$ converge to $T$ such that $L_{t_{i_1}}$ and $L_{t_{i_2}}$ converge to the varifolds $L_1$ and $L_2$. Let us denote the radon measures associated to the varifolds by $\mu_1$ and $\mu_2$ respectively. 
    \begin{claim}
        $\mu_1 = \mu_2$. 
    \end{claim}
    \begin{claimproof}
    Let $f$ be any compactly supported smooth function, $f \geq 0.$ By [\citealp{ilmanen1994elliptic}, 7.2], there exists $C_f$ such that $\int_{L_t} f - C_f t$ is decreasing in $t$. Hence, for any $t_{j_2},$
    \begin{align*}
        \mu_1(f) &= \lim\limits_{t_{i_1}} \Big(\int_{L_{t_{i_1}}} f - C_f t_{i_1} \Big) + C_f T \\
        &\leq \int_{L_{t_{j_2}}} f - C_f t_{j_2} + C_f T \\
        &= \int_{L_{t_{j_2}}} f + C_f (T - t_{j_2}).
    \end{align*}
    Taking $j_2 \to \infty$ to both sides we have 
    $$\mu_1(f) \leq \lim\limits_j  \int_{L_{t_{j_2}}} f  =  \mu_2(f).$$
    Now let $\epsilon >0.$ Then there exists a $t_{j_1}$ such that $T - t_{j_1} < \epsilon$ and $|\mu_1(f) - \int_{L_{t_{j_1}}} f| < \epsilon.$ Choose any $t_{i_2} \in [t_{j_1}, T).$ Then
    \begin{align*}
        \int_{L_{t_{i_2}}} f - \mu_1(f) &\leq \int_{L_{t_{i_2}}} f  - \int_{L_{t_{j_1}}} f + \epsilon \\
        &\leq \lim\limits_{t \to {t^{+}_{i_2}}} \Big(\int_{L_t} f - C_f t \Big) + C_f t_{i_2} - \int_{L_{t_{j_1}}} f + \epsilon \\
        &\leq \int_{L_{t_{j_1}}} f - C_f t_{j_1} + C_f t_{i_2} - \int_{L_{t_{j_1}}} f  + \epsilon \\
        &= C_f (t_{i_2} - t_{j_1}) + \epsilon.
    \end{align*}
    where in the second line we have used [\citealp{ilmanen1994elliptic}, $7.2$]. Now take $i_2 \to \infty$ and then $\epsilon \to 0$ to get 
    $$\mu_2(f) \leq \mu_1(f).$$
    Combining both we can conclude that $\mu_1(f) = \mu_2(f).$ As $f$ arbitrary, so the claim follows.
    \end{claimproof} \\
    So for any sequence $t_i \to T,$ there exists a subsequence of $t_i$ which converges to the same varifold. This proves that $L_t$ converges to a unique varifold as $t \to T$.
    
\end{proof}

Now pick a sequence $t_i \to T$ and extract a subsequence (still denoted by index $i$) to get a limiting current at time $T$, we will call it $\widehat{L}_T$, which apparently depends on the the sequence $t_i \to T.$ We will denote the radon measure associated to the varifold $L_T$ and the current $\widehat{L}_T$ by $\mu_{L_T}$ and $\mu_{\widehat{L}_T}$ respectively. Note that the measures $\mu_{L_T}$ and $\mu_{\widehat{L}_T}$ might not be the same. We will be more interested in the limiting current as for the later part we need to deal with currents instead of the varifolds.

\begin{thm} {\label{thm 3.2}}
    Let $(L_t)_{t \in [0,T)}$ be a solution to the Lagrangian mean curvature flow with uniformly bounded mean curvature and for every sufficiently small $r > 0,$ if there exists a $t_r < T$ close to $T$, depending on $r$ such that for $t_r \leq t < T$, there is only one connected component of $L_t \cap B_{4r}(0)$ that intersects $B_r(0)$. Then there is a one to one correspondence between the tangent flows at the space-time $(x_0,T)$ and the static flow of the tangent cone of the varifold associated to the current $\widehat{L}_T$ at $x_0$. 
\end{thm}

We are going to prove Theorem $\ref{thm 3.2}$ in the end of this section. First we will prove that the tangent flows of the flow are the same as the static flows of tangent cones of the varifold $L_T$. Then we will show the tangent cones of the varifold $L_T$ are the tangent cones of the varifold associated to the current $\widehat{L}_T$. Note that this is not an immediate consequence because those varifolds might not be same.

Without loss of generality, we can take $x_o$ to the origin. Fix $r>0$ small. As $L_t$ develops a singularity at the origin, then $\exists \;t_0$ such that if $t \geq t_0$ then $L_t \cap B_r(0)$ is non-empty.

We will first focus on the number of components of $L_t \cap B_{4r}(0)$ that intersects $B_r(0).$ Pick $t_1$ such that $t_0 \leq t_1 <T.$ Either $L_{t_1} \cap B_{r}(0)$ is connected in $B_{4r}(0)$ or there are at least two connected components of $L_{t_1} \cap B_{4r}(0)$ that intersects $B_r(0).$ Let $C_1^{t_1}$ and $C_2^{t_1}$ be two such connected components. 

Now we will show if $t_1$ is sufficiently close to $T,$ then $C_1^{t_1} \cap B_r(0)$ and $C_2^{t_1} \cap B_r(0)$ never intersect under the flow at a later time for $t_1 \leq t < T.$ Choose $t_1 \geq \text{max} \{t_0, T - \frac{r}{4\Lambda} \}$. If we denote $F$ be the embedding map, then we have, 
$$-2 \Lambda \leq \frac{\partial}{\partial t}|F(p,t)| \leq 2 \Lambda.$$
So it follows that, for any $t_1,t_2$
$$|F(p,t_1)| \leq |F(p,t_2)| + 2 \Lambda |t_1-t_2|.$$
Hence if we denote $C_x(B_r(x) \cap L)$ by the connected component of $B_r(x) \cap L$ containing $x$, then we have
$$C_{F(p,t_2)} \big(B_r(F(p,t_2)) \cap L_{t_2}\big) \subset C_{F(p,t_1)} \big(B_{3r/2}(F(p,t_1)) \cap L_{t_1}\big). $$
Now if $F(p,t_1) \in C_1^{t_1} \cap B_r(0)$ and $F(q,t_1) \in C_2^{t_1} \cap B_r(0)$ then this implies
$$F(p,t_2) \in C_{F(p,t_1)} \big(B_{3r/2}(F(p,t_1)) \cap L_{t_1}\big),$$ 
but 
$$C_{F(p,t_1)} \big(B_{3r/2}(F(p,t_1)) \cap L_{t_1}\big) \subset C_1^{t_1} \cap B_{4r}(0)$$
therefore,
$$F(p,t_2) \in C_1^{t_1} \cap B_{4r}(0).$$
By a similar argument,
$$F(q,t_2) \in C_2^{t_1} \cap B_{4r}(0).$$
So, $C_1^{t_1} \cap B_r(0)$ and $C_2^{t_1} \cap B_r(0)$ never intersect under the flow at a later time for $t_1 \leq t < T.$
Also, note that we only care about what happens in $B_r(0)$ at time $t_1$ because, if at time $t_1, F(p,t_1) \notin B_r(0),$ then for all later time $t_2,$
$$|F(p,t_2)| \geq r - \frac{r}{2} = \frac{r}{2}$$

So, from the above discussion we can conclude that if the mean curvature stays uniformly bounded, then for every small $r >0$ there exists a $t_r < T$ depending on $r$ such that for $t_r \leq t < T$, there is only one connected component of $L_t \cap B_{4r}(0)$ that intersects $B_r(0)$ or $L_{t_r} \cap B_r (0)$ has more than one connected component in $B_{4r}(0)$ and they remain disconnected in $B_r(0)$ for all $t_r \leq t < T$.

From now on, we will assume that there is only one connected component of $L_t \cap B_{4r}(0)$ that intersects $B_r(0)$. First we will show that in this case, $L_t$ satisfies an almost minimizing property near the origin in the following sense,

\begin{lemma} {\label{lemma 3.3}}
    We have the following,
    $$M_{B_r(0)} (L_t) \leq \frac{1}{1-Cr^2} M_{B_r(0)} (L_t + \partial S)$$ 
    for any $3$ dimensional current $S$ supported in $B_r(0).$ 
\end{lemma}

\begin{proof}
    First note that, we have a uniform lower bound 
    $$\mathcal{H}^2(\hat{B}_r(x)) \geq Kr^2$$ 
    for the intrinsic ball of radius $r$ in $L_t$ centered at any $x \in L_t \cap B_r(0).$ This has been proved in [\citealp{neves2007singularities}, Lemma $7.2$]. Also, by bounded area ratio, we  have an upper bound for the area, 
    $$\mathcal{H}^2(L_t \cap B_r(0)) \leq Cr^2.$$ 
    So, if we take any two points $x_1, x_2 \in L_t \cap B_r(0),$ then the points can be joined by a curve in $L_t$ has length at most $Cr$.\\
    Fix a point $x \in L_t \cap B_r(0)$. We can assume, $\theta_{L_t}(x) = 0.$ Now take any point $y \in L_t \cap B_r(0).$ As $|H(L_t)| \leq \Lambda $ and $H = J\nabla \theta,$ so 
    $$|\theta_t(x) - \theta_t(y)| \leq C d_{L_t}(x,y).$$
    As, 
    $\Omega|_{L_t} = e^{i\theta} Vol_{L_t},$ so 
    $$Re \Omega|_{L_t} = (cos{\theta}) Vol_{L_t}.$$
    Now, 
    $$cos (\theta) \geq 1 -\frac{\theta^2}{2} = 1 - Cr^2,$$
    so on $B_r(0),$ we have 
    $$Vol_{L_t} \leq \frac{1}{1- Cr^2}Re \Omega|_{L_t}.$$ 
    Therefore,
    \begin{align*}
        \mathcal{H}^2 (L_t \cap B_r(0)) &= \int_{B_r(0)} Vol_{L_t} \\
        &\leq \frac{1}{1-Cr^2}\int_{B_r(0)} Re \Omega|_{L_t} \\
        &= \frac{1}{1-Cr^2} \int_{(L_t + \partial S)} 1_{B_r(0)} Re \Omega - \frac{1}{1-Cr^2} \int_{\partial S} 1_{B_r(0)} Re \Omega \\
        &= \frac{1}{1-Cr^2} \int_{(L_t + \partial S)} 1_{B_r(0)} Re \Omega - \frac{1}{1-Cr^2}  S (d Re \Omega|_{B_r(0)}) \\
        &\leq \frac{1}{1-Cr^2} \mathcal{H}^2((L_t + \partial S)(B_r(0)). 
    \end{align*}
    
\end{proof}

Now using Lemma \ref{lemma 3.3}, we will show the almost minimizing condition holds for the limiting current $\widehat{L}_T$ as well. Note that, this is not a direct consequence of convergence of currents as convergence of currents only preserve the lower semi-continuity of mass. 

\begin{lemma} {\label{lemma 3.4}}
    $\widehat{L}_T$ is almost minimizing near the origin as in Lemma \ref{lemma 3.3} i.e.
    $$M_{B_r(0)} (\widehat{L}_T) \leq \frac{1}{1-Cr^2} M_{B_r(0)} (\widehat{L}_T + \partial S)$$
    for any $3$ dimensional current $S$.
\end{lemma}

\begin{proof}
    As $L_{t_i}$ converges $\widehat{L}_T$ in the sense of current, so there are currents $A_i$ and $B_i$ such that $M (A_i) + M (B_i) \rightarrow 0$
    $$(\widehat{L}_T - L_{t_i}) \cap B_r(0) = \partial A_i + B_i.$$
    Let $S$ be any three dimensional integer multiplicity current supported in $B_r(0)$.
    From Lemma \ref{lemma 3.3} we have,    
    $$M_{B_r(0)} (L_{t_i}) \leq \frac{1}{1-Cr^2} M_{B_r(0)} (L_{t_i} + \partial S).$$
    Hence,
    \begin{align*}
        M_{B_r(0)} (\widehat{L}_T + \partial S) &= M_{B_r(0)} (L_{t_i} + \partial A_i + B_i + \partial S) \\
        &\geq M_{B_r(0)} (L_{t_i} + \partial(A_i + S)) - M_{B_r(0)} (B_i) \\
        &\geq (1-Cr^2) M_{B_r(0)} (L_{t_i}) - M_{B_r(0)} (B_i).
    \end{align*}
    So taking $\liminf$ to the both sides, we get 
    $$M_{B_r(0)} (\widehat{L}_T + \partial S) \geq (1-Cr^2) \liminf M_{B_r(0)} (L_{t_i}).$$
    By the lower semi-continuity of mass with respect to weak convergence we have
    $$M_{B_r(0)} (\widehat{L}_T) \leq \liminf M_{B_r(0)} (L_{t_i}).$$
    Combining both inequalities give,
    $$M_{B_r(0)} (\widehat{L}_T) \leq \frac{1}{1-Cr^2} M_{B_r(0)} (\widehat{L}_T + \partial S).$$
\end{proof}

\begin{lemma} {\label{lemma 3.5}}
    We have the following,
    $$\limsup M_{B_r(0)} (L_{t_i}) \leq \frac{1}{1-Cr^2} M_{B_r(0)} (\widehat{L}_T).$$ 
\end{lemma}

\begin{proof}
    The proof is exactly like before. In the last proof, taking $S= 0,$ we get
    \begin{align*}
        M_{B_r(0)} (\widehat{L}_T) &=  M_{B_r(0)} (L_{t_i} + \partial A_i + B_i) \\
        &\geq  M_{B_r(0)} (L_{t_i} + \partial A_i) -  M_{B_r(0)} (B_i) \\
        &\geq (1-Cr^2)  M_{B_r(0)} (L_{t_i}) -  M_{B_r(0)} (B_i). 
    \end{align*}
    Now take $\limsup$ to the both sides to conclude the lemma. \\
\end{proof}
    
Now we will show that any tangent flow is a static flow if the mean curvature remains uniformly bounded, which is our next lemma. 

\begin{lemma} {\label{lemma 3.6}}
    Any tangent flow is a static flow if the mean curvature remains uniformly bounded throughout the flow.
\end{lemma}

\begin{proof}
    Consider a cut-off function $\chi$ supported in the ball of Radius $R$. Let us denote $L^i_s$ by  
    $$L^i_s := \lambda_i \big(L_{T+ \frac{s}{\lambda_i^2}} -x_0 \big).$$
    For $-\lambda_i^2T \leq s < 0,$ we have
    \begin{align*}
        \Big| \partial_s \int_{L^i_s} \chi \Big| &= \Big| \int_{L^i_s} \big(\chi|H|^2 +  \langle D\chi,H \rangle \big)\Big| \\
        &\leq C \int_{L^i_s} \big( |H|^2 + |H| \big).
    \end{align*}
    So for any $-\lambda_i^2T \leq s_1 < s_2 < 0$,  
    \begin{align*}
        \Big|\int_{L^i_{s_1}} \chi - \int_{L^i_{s_2}} \chi \Big| \leq C \int_{s_1}^{s_2} \int_{L^i_t} \big( |H|^2 + |H| \big).
    \end{align*}
     As $|H(L_t)| \leq \Lambda$, so $|H(L^i_{t})| \leq \frac{\Lambda}{\lambda_i}$. So if $\lambda_i$ goes to $\infty,$ the right hand side goes to $0$. Hence, $L^i_s$ converges to a static flow.
\end{proof}
Let us now define the varifolds $L_0^i$ by 
    $$L_0^i := \lambda_i (L_T - x_0)$$
Now for any sequence $i \to \infty, L^i_s$ has a subsequence $L^{i_k}$ converges weakly to a Brakke flow. By compactness theorem of Radon measure, there exists a subsequence $i_{k_j}$ such that $\mu_{L^{i_{k_j}}_T}$ converges to a Radon measure $\mu.$ We will denote the  subsequence $i_{k_j}$ by $i.$ The next lemma says the Brakke flow is the static flow $L_s$ and at each time slice, $L_s$ is the same stationary varifold $\mu.$ 

\begin{lemma} {\label{lemma 3.7}}
The tangent flow of $L_t$ and the static flow generated by the tangent cones of $L_T$ are same.
\end{lemma}
         
\begin{proof}
        Let $\epsilon >0.$ Note that for any $i,$ by Varifold convergence we have there exists $s_0$ close to $0$ we have
        $$\Big|\int_{L^i_{s_0}} \chi - \int_{L^i_0} \chi \Big| < \frac{\epsilon}{2}.$$
        We have shown, for large $i$
        $$\Big|\int_{L^i_{s}} \chi - \int_{L^i_{s_0}} \chi \Big| \leq \frac{\epsilon}{2}.$$
        Combining both inequalities, we have for any $s,$
        $$\Big|\int_{L^i_s} \chi - \int_{L^i_0} \chi \Big| \leq \epsilon,$$
        which is what we wanted to prove.
\end{proof} 
Let us now similarly define the currents 
$$\widehat{L}_0^i := \lambda_i (\widehat{L}_T - x_0).$$
By compactness theorem of Radon measure, we know for any sequence $i \to \infty$, there exists a subsequence $i_k$ such that $\mu_{L^{i_k}_T}$ converges to a Radon measure $\mu.$ Now we will choose a subsequence $i_{k_j}$ of $i_k$ such that $\mu_{\widehat{L}^{i_{k_j}}_T}$ converges to a Radon measure $\widehat{\mu}.$ We will just denote the subsequence by $i_{k_j}$ by $i$ as before. Next we will show that 
$$\mu = \widehat{\mu},$$ 
which is our next claim. 

\begin{lemma} {\label{lemma 3.8}}
There is a one to one correspondence between the tangent cones of $L_T$ and the tangent cones of the varifold associated to current $\widehat{L}_T.$
\end{lemma}

\begin{proof}
As $L_t$ converges to $\widehat{L}_T$ in the sense of current, so by weak convergence 
$$\mu_{\widehat{L}_0^i} \big(B_r(0) \big) \leq \liminf \;\mu_{L_t^i} \big(B_r(0) \big).$$ 
Hence, $\mu_{\widehat{L}_0^i} \big(B_r(0)\big) \leq \mu_{L_T^i} \big(B_r(0)\big)$ which implies 
$\widehat{\mu} (B_r(0)) \leq \mu(B_r(0))$. Now, using Lemma $\ref{lemma 3.5}$, 
\begin{align*}
         \limsup \;\mu_{L_t^i}(B_r(0)) &= \limsup\;\mu_{\lambda_iL_{T+ \frac{t}{\lambda_i^2}}}(B_r(0)) \\
         &= \limsup\;\lambda_i^n \mu_{L_{T+ \frac{t}{\lambda_i^2}}}(B_{\frac{r}{\lambda_i}}(0)) \\
         &\leq \frac{\lambda_i^n}{1-\frac{Cr^2}{\lambda_i^2}} \mu_{\widehat{L}_T}(B_{\frac{r}{\lambda_i}}(0)) \\
         &\leq \frac{1}{1-\frac{Cr^2}{\sigma_i^2}} \mu_{\widehat{L}_0^i}(B_r(0)). 
\end{align*}  
So this implies $\mu \big(B_r(0) \big) \leq \widehat{\mu} \big(B_r(0) \big).$ Hence, 
$$\mu \big(B_r(0) \big) = \widehat{\mu} \big(B_r(0) \big).$$ 
But as ${\mu}_{\widehat{L}_T}(A) \leq \mu_{L_T}(A)$ for all open sets $A$, so $\widehat{\mu}(A) \leq \mu(A)$, hence by Radon-Nikodym theorem,
$$\frac{d\widehat{\mu}}{d\mu} \leq 1$$ 
$\mu$ almost everywhere, but it can not be strictly less than $1$ on a set of positive measure as otherwise it will contradict $\mu\big(B_r(0)\big) = \widehat{\mu} \big(B_r(0)\big).$ So we can conclude that $\widehat{\mu} = \mu.$ 
\end{proof} 
    
We have shown in Lemma $\ref{lemma 3.7}$ the correspondence between the tangent flows and the static flow of the tangent cones of the limiting varifold $L_T$ and in Lemma $\ref{lemma 3.8}$, the correspondence between the the static flow of the tangent cones of the limiting varifold and the static flow of the tangent cones of the varifold associated to limiting current $\widehat{L}_T$. As the varifold does not depend on the sequence, so for any such limiting current the static flows of the tangent cones of the varifold associated to it are same even though the limiting current might be different. Combining these two gives us a one to one correspondence between the tangent flows and the static flows of the tangent cones of varifold associated to a limiting current. This proves Theorem $\ref{thm 3.2}.$

\section{Epiperimetric type inequality}\label{}

We will denote $x^T$ for the projection of the position vector $x$ on the approximate tangent space (for definition see \cite{simon1983lectures}) to $\widehat{L}_T$ at $x$. For $r$ small enough, $$\frac{1}{1-Cr^2} \leq (1+2Cr^2).$$
From now on, we will denote $1+2Cr^2$ by $C_r$.

In this section we will prove that limiting current $\widehat{L}_T$ has a unique tangent cone, which is our next theorem. Next we will show that this is equivalent to prove our main Theorem \ref{thm 1.1}.

\begin{thm} {\label{thm 4.1}}
    There is an area minimizing cone such that $r^{-1} \widehat{L}_T$ converges to the cone as $r \to 0$ in the sense of current.
\end{thm}

This theorem will be proved in the end of this paper. The idea of the proof is to use epiperimetric inequality similar to White \cite{white1983tangent} (also see De Lellis-Spadaro-Spolaor \cite{de2017uniqueness}). First we need a monotonicity formula for $\widehat{L}_T$. This is similar to the the monotonicity formula for Minimal surfaces (also see \cite{de2017uniqueness}, Proposition $2.1$).

\begin{prop} {\label{prop 4.2}}
    There exist $r > 0,$ so that if $0< t<s<r,$ then 
    $$\frac{M_{B_s(0)} (\widehat{L}_T)}{\omega_2 s^2} - \frac{M_{B_t(0)} (\widehat{L}_T)}{\omega_2 t^2} + \epsilon(s,t) \geq \frac{1}{3\omega_2}\int_{B_s(0) \setminus B_t(0)} \frac{|x^\perp|^2}{|x|^4} \;d \mu_{\widehat{L}_T} $$
    where,
    $$\epsilon(s,t) = \frac{1}{3\omega_2} \int_{B_s(0) \setminus B_t(0)} \frac{(C_r^2 - 1)|x^T|^2}{|x|^4} \;d \mu_{\widehat{L}_T} .$$
\end{prop}

\begin{proof}
    By Lemma \ref{lemma 3.4}, we have 
    $$M_{B_r(0)} (\widehat{L}_T) \leq \frac{r}{2} C_r M(\partial (\widehat{L}_T \mres B_r(0)).$$
    Let 
    $$f(r) := M_{B_r(0)} (\widehat{L}_T).$$ 
    As $f$ increasing, so $Df$ is non-negative measure. Let $f'$ be the absolutely continuous part and $\mu_s$ be the singular part of one dimensional Lebesgue measure with respect to $Df.$ Hence, 
    $$\frac{\mu_s}{r^2} + \frac{f'(r)}{r^2} - \frac{1}{r^2}C_r M (\partial (\widehat{L}_T \mres B_r(0)) \leq \frac{Df}{r^2} - \frac{2}{r^3} M_{B_r(0)} (\widehat{L}_T).$$
    So,
    $$\int_t^s \frac{1}{r^2} d\mu_s + \int_t^s \bigg(\frac{f'(r)}{r^2} - \frac{1}{r^2}C_r M (\partial (\widehat{L}_T \mres B_r(0))\bigg) dr \leq \frac{f(s)}{s^2} - \frac{f(t)}{t^2}.$$
    Now,
    \begin{align*}
        \int_t^s \bigg(\frac{f'(r)}{r^2} - \frac{1}{r^2}C_r M (\partial (\widehat{L}_T \mres B_r(0))\bigg) dr &= \int_t^s \frac{1}{r^2} \int \frac{|x| - C_r|x^T|}{|x^T|} \;d \mu_{\widehat{L}_T} \\
        &= \int_t^s \frac{1}{r^2} \int \frac{|x|^2 - C_r^2|x^T|^2}{|x^T|\big(|x| + C_r|x^T|\big)} \;d \mu_{\widehat{L}_T} . 
    \end{align*}
    Choose $r$ small enough so that $C_r \leq 2,$ then
    $$|x| + C_r|x^T| \leq (1+C_r)|x| \leq 3|x|.$$
    So we have,
    \begin{align*}
        \int_t^s \frac{1}{r^2} \int \frac{|x|^2 - C_r^2|x^T|^2}{|x^T|\big(|x| + C_r|x^T|\big)} \;d \mu_{\widehat{L}_T} &= \int_t^s \frac{1}{3r^4} \int \frac{\big(|x|^2 - C_r^2|x^T|^2\big)|x|}{|x^T| } \;d \mu_{\widehat{L}_T} \\
        &= \int_{(B_s(0) \setminus B_t(0)) \cap \{|x^T| >0\}} \frac{|x|^2 - C_r^2|x^T|^2}{3|x|^4} \;d \mu_{\widehat{L}_T} .
    \end{align*}
    Now,
    \begin{align*}
        \int_t^s \frac{1}{r^2} d\mu_s &= \lim\limits_{N \to \infty} \sum_{i=1}^N \frac{1}{(t+\frac{i-1}{N}(s-t))^2} \int^{t+\frac{i}{N}(s-t)}_{t+\frac{i-1}{N}(s-t)}
        d\mu_s \\
        &\geq \lim\limits_{N \to \infty} \sum_{i=1}^N \frac{1}{(t+\frac{i-1}{N}(s-t))^2} \int_{(B_{t+\frac{i}{N}(s-t)}(0) \setminus B_{t+\frac{i-1}{N}(s-t)}(0)) \cap \{x^T = 0\}} d\mu_{\widehat{L}_T}\\
        &\geq \int_{(B_s(0) \setminus B_t(0)) \cap \{x^T = 0\}} \frac{1}{|x|^2} \; d\mu_{\widehat{L}_T} \\
        &= \int_{(B_s(0) \setminus B_t(0)) \cap \{x^T = 0\}} \frac{|x^\perp|^2}{|x|^4} \;d\mu_{\widehat{L}_T}\\
        &\geq \int_{(B_s(0) \setminus B_t(0)) \cap \{x^T = 0\}} \frac{|x|^2 - |x^T|^2}{3|x|^4} \;d\mu_{\widehat{L}_T}.
    \end{align*}
    So adding these two inequalities we have our result.
\end{proof}

\begin{remark} \label{rmk 4.2}
    Note that,
    \begin{align*}
    \epsilon(s,t) &= \frac{1}{3\omega_2} \int_{B_s(0) \setminus B_t(0)} \frac{(C_r^2 - 1)|x^T|^2}{|x|^4} \;d \mu_{\widehat{L}_T} \\
    &\leq C \int_{B_s(0) \setminus B_t(0)} \frac{Cr^2}{|x|^2} \;d \mu_{\widehat{L}_T} \\
    &\leq Cs^2,
\end{align*}
where $C$ is a constant which varies line to line.
\end{remark}

Now we will prove the epiperimetric inequality for limiting current $\widehat{L}_T$ in the next lemma. To prove this we need the main result from \cite{white1983tangent} (also see \cite{de2017uniqueness}, Lemma $3.3$). The following result and the proof are very similar to Proposition $3.4$ in \cite{de2017uniqueness}. The only difference is that in \cite{de2017uniqueness}, they assumed the similar almost minimizing property for any point in support of $M$ which was needed for the proof. But we have the almost minimizing property just near the origin. We will overcome this difficulty with an extra assumption that is the condition $2$ in the lemma. Later, we will see current $\widehat{L}_T$ satisfies all conditions as in the lemma.

\begin{lemma} \label{lemma 4.4}
    Let $T$ be an area minimizing cone and $r >0$. Then there exists an $\epsilon > 0$ depending on $T$ such that if $M$ satisfies the following conditions,
    \begin{enumerate} [before=\normalfont\normalsize]
        \item it has a uniform mass bound for any compact set, 
        \item there exists a $C>0$ such that for any $\delta > 0,$ there exists $r_\delta >0$ depending on $\delta$ such that  
        $$\mu_{B_\delta (x)} (r^{-1}M) \geq C\delta^2$$
        for every $0< r \leq r_\delta$ and for every $x \in$ spt $r^{-1}M \cap B_2(0),$
        \item it has an almost minimizing property near the origin in the sense of Lemma $\ref{lemma 3.3}$,
        \item the flat norm of the current $\big(r^{-1}M - T\big) \;\mres B_2(0)$ is smaller than $\epsilon >0$.
    \end{enumerate}
    Then we have the following, 
    \begin{align*}
        M_{B_1(0)}(r^{-1}M) - M_{B_1(0)}(T) &\leq C_r(1- \epsilon)\Big(M_{B_1(0)} \big(0 \; {\RangeX}\;\partial (r^{-1}M \cap B_1(0))\big) \;-\\ 
        &\hspace{0.5 cm} M_{B_1(0)}(T) \Big) +  (C_r-1) M_{B_1(0)}(T).
    \end{align*}
\end{lemma}

\begin{proof}
    The proof is by contradiction. Let there exists a sequence $r_n \to 0$ such that the flat norm of the currents $\big(r_n^{-1}M_n - T\big) \cap B_2(0)$ are smaller than $\frac{1}{n}$ but 
    \begin{align*}
        M_{B_1(0)}(r_n^{-1}M_n) - M_{B_1(0)}(T) &\geq C_{r_n} \bigg(1- \frac{1}{n} \bigg) \Big(M_{B_1(0)} \big(0 \; {\RangeX}\;\partial (r_n^{-1}M_n \cap B_1(0))\big) \;-\\ 
        &\hspace{0.5 cm} M_{B_1(0)}(T) \Big) +  (C_{r_n}-1) M_{B_1(0)}(T).
    \end{align*}
    As, $M_{B_{2r_n}(0)}(M_n) \leq 4Cr_n^2,$ so we have uniform bound for $B_2(0)$ of $r_n^{-1}M_n$. By slicing theorem (see \cite{simon1983lectures}) and the fact  $r_n^{-1}M_n$ converges to $T$ in the sense of flat norm on $B_2(0)$, we can conclude that $\exists \;r_0 >1$ such that the flat norm of the current $(r_n^{-1}M_n - T) \mres B_{r_0}(0)$ goes to 0. So there are integral currents $A_n$ and $B_n$ such that
    $$(r_n^{-1}M_n - T) \cap B_{r_0}(0)  = \partial A_n + B_n$$
    $$(M_n - r_nT) \cap B_{r_nr_0}(0) = \partial(r_nA_n) + r_nB_n.$$
    So by the almost minimizing property of $M_n$,
    \begin{align*}
        M_{B_r(0)} (r_n^{-1}M_n)  &= \frac{1}{r_n^2} M_{B_{r_nr_0}(0)} (M_n) \\
        &\leq \frac{C_{r_nr_0}}{r_n^2} M_{B_{r_nr_0}(0)} \big(M_n + \partial(-r_nA_n)\big)\\
        &=  \frac{C_{r_nr_0}}{r_n^2} M_{B_{r_nr_0}(0)} \big(r_n(T + B_n) \big) \\
        &\leq C_{r_nr_0} M_{B_{r_0}(0)} T+B_n. 
    \end{align*} 
    Hence,
    $$M_{B_{r_0}(0)} (r_n^{-1}M_n) \leq C_{r_nr_0} M_{B_{r_0}(0)} (T) + C_{r_nr_0}  M(B_n).$$
    When $n \to \infty, C_{r_nr_0} \to 1$ and $M(B_n) \to 0,$ so
    $$\limsup M_{B_{r_0}(0)} (r_n^{-1}M_n)  \leq M_{B_{r_0}(0)} (T).$$
    By weak convergence,
    $$M_{B_{r_0}(0)} (T) \leq \limsup  M_{B_{r_0}(0)} (r_n^{-1}M_n).$$
    So we can conclude that 
    $$M_{B_{r_0}(0)} (r_n^{-1}M_n) \to M_{B_{r_0}(0)} (T).$$
    Now we can take a subsequence of $r_n^{-1}M_n$ which converges as a radon measure because of condition $1$. Now as $M_{B_{r_0}(0)} (r_n^{-1}M_n) \to M_{B_{r_0}(0)} (T),$ so we can conclude that $r_n^{-1}M_n$ (taking a subsequence and reindexing it to the original sequence) converges to $T$ as a radon measure on $B_2(0).$ The proof follows exactly as in the proof of the claim of our paper. Now $T$ is a minimal cone in $\mathbb{C}^2$, so the measure of $T$ on the boundary of the unit ball vanishes, hence we can conclude
    $$M_{B_1(0)}(r_n^{-1}M_n ) \to M_{B_1(0)}(T).$$ 
    The next step is to show that 
    $$M \big(\partial(r_n^{-1}M_n \cap B_1(0)) \big) \to M \big (\partial(T \cap B_1(0)) \big).$$
    As, $M_{B_1(0)}(r_n^{-1}M_n) - M_{B_1(0)}(T) \to 0,$ so by our assumption
    $$M \Big(0 \; {\RangeX}\;\partial (r_n^{-1}M_n \cap B_1(0))(B_1(0)) - T(B_1(0))\Big)\to 0.$$ 
    As $T$ is area minimizing cone, so 
    $$T \mres B_1(0) = 0 \; {\RangeX}\;\partial (T \cap B_1(0)),$$
    and therefore, 
    $$M \big(\partial(r_n^{-1}M_n \cap B_1(0)) \big) \to M \big(\partial(T \cap B_1(0)) \big).$$
    Now, as the flat norm of the boundary of the current $(r_n^{-1}M_n - T) \cap B_1(0)$ is bounded by the mass
    $$M \Big(0 \; {\RangeX}\;\partial (r_n^{-1}M_n \cap B_1(0))(B_1(0)) - T(B_1(0))\Big),$$ 
    so we can conclude that the flat norm of the boundary of the current $(r_n^{-1}M_n - T) \cap B_1(0)$ converges to $0$ as well. 
    
    Now we will show that the spt $r_n^{-1}M_n \cap B_2(0)$ converges to spt $T \cap B_2(0)$ in the Hausdroff sense.
    
    Let $x \in spt\; T \cap \overline{B_1(0)}$. Choose any $\delta > 0.$ For any set $A,$ define 
    $$B_\delta (A) = \{x: d(x,A) < \delta\}.$$
    As $M_{B_\delta(x)} (T) > 0$, so for large $n, M_{B_\delta(x)} (r_n^{-1}M_n) > 0$ as $r_n^{-1}M_n$ converges to $T$ in measure on $B_{r_0}(0)$ for some $r_0 >1$. Hence, there exists $x_n \in B_\delta(x)$ such that $x_n \in spt \;r_n^{-1}M_n.$ Therefore, for large $n$
    $$spt \;T \subset B_\delta (spt \;r_n^{-1}M_n).$$
    Now suppose there exists a sequence $x_n \in\; \overline{B_1(0)} \cap spt \;r_n^{-1}M_n \setminus B_\delta(spt \;T).$ As $x_n \in\; spt \;r_n^{-1}M_n,$ so by condition $2$, there exists a $N \in \mathbb{N}$ such that if $n \geq N$, then $$M_{B_{\delta/2}(x_n)}(r_n^{-1}M_n) \geq K\delta^2$$ 
    for some $K$. Now we can find a subsequence of $x_n$ which coverges to $x \in \overline{B_1(0)}$, as $x_n \notin B_\epsilon (spt T)$ so $M_{B_\delta(x)}(T) = 0$ which implies 
    $$M_{B_\delta(x_n)}(r_n^{-1}M_n) \to 0.$$
    But, for large $n, B_{\delta/2}(x_n) \in B_{\delta}(x)$. Therefore, for large $n,$ 
    $$M_{B_\delta(x)} (r_n^{-1}M_n) \geq K\delta^2,$$ 
    which is a contradiction.
    
    Hence,
    $$spt \;r_n^{-1}M_n \subset B_\delta (spt \;T).$$
    So from Lemma $3.3$ in \cite{de2017uniqueness} (See also \cite{white1983tangent}), we can conclude that there exist integral currents $H_n$ with $(H_n - r_n^{-1}M_n) \cap B_1(0)$ has zero boundary such that, 
    $$M_{B_1(0)}(H_n) - M_{B_1(0)}(T) \leq (1-\epsilon)\Big(M_{B_1(0)} \big(0 \; {\RangeX}\;\partial (r^{-1}M \cap B_1(0))\big) - M_{B_1(0)}(T) \Big).$$
    But, by the almost minimizing property, 
    $$M_{B_1(0)} (r_n^{-1}M_n) \leq C_{r_n} M_{B_1(0)} (H_n).$$
    So we have,
    \begin{align*}
        M_{B_1(0)} (r_n^{-1}M_n) - M_{B_1(0)}(T) &\leq C_{r_n} M_{B_1(0)} (H_n) - M_{B_1(0)} (T) \\
        &= C_{r_n} M_{B_1(0)} (H_n) - C_{r_n} M_{B_1(0)}(T) + (C_{r_n} -1) \\
        &\hspace{0.5 cm} M_{B_1(0)}(T) \\
        &\leq C_{r_n} (1- \epsilon) \big(0\; {\RangeX}\; M_{B_1(0)} \big(\partial (r_n^{-1} M \cap B_1(0)) \;-\\
        &\hspace{0.4 cm} M_{B_1(0)}(T) \big)+(C_{r_n} -1) M_{B_1(0)}(T).
    \end{align*}
    Now choose $n$ large so that $\frac{1}{n} < \epsilon,$ then we will get a contradiction.
\end{proof}

Now we will use Lemma \ref{lemma 4.4} to prove our next proposition which establishes a one to one correspondence between the tangent flow of the original flow and the tangent cone of limiting current. 

\begin{prop} {\label{prop 4.5}}
    If $(L_t)_{t \in [0,T)}$ be a solution to the Lagrangian mean curvature flow with uniformly bounded mean curvature in $\mathbb{C}^2$ and for every sufficiently small $r > 0,$ if there exists a $t_r < T$ close to $T$, depending on $r$ such that for $t_r \leq t < T$, there is only one connected component of $L_t \cap B_{4r}(0)$ that intersects $B_r(0)$, then there is a one to one correspondence between the tangent flows at $(x_0,T)$ and the tangent cone of current $\widehat{L}_T$ at $x_0$. 
\end{prop}

\begin{proof}
For any sequence $i \to \infty$, as before we can extract a subsequence $i_k$ such that $L_s^{i_k}$ converges weakly to a Brakke flow. By compactness theorem of integral current, there exists a subsequence $i_{k_j}$ such that $\widehat{L}_0^{i_{k_j}}$ converges to a current $\widehat{T}$ and by compactness theorem of the integral varifold, the varifold associated to the current $\widehat{L}_0^{i_{k_j}}$ converges to a varifold $T$. We will denoted the subsequence by $i.$ We already have shown in Theorem $\ref{thm 3.2}$ that the Brakke flow is static flow, and at each time slice it is the same stationary varifold $T$. So all we need to show is that the varifold associated to the current $\widehat{T}$ is the same as the varifold $T.$ 

First we will show that $\widehat{L}_T$ satisfies all the conditions given in Lemma $\ref{lemma 4.4}.$

The current $\widehat{L}_T$ has uniform mass bound for any compact set $K$ because $L_t$ has uniformly bounded area ratios and 
$$M_K(\widehat{L}_T) \leq \liminf M_K(L_t).$$

We already have shown in Lemma \ref{lemma 3.4} that $\widehat{L}_T$ satisfies the almost minimizing condition. 

As $\widehat{L}^i_0$ converges to $T$ as a current, so the flat norm of $\widehat{L}^i_0 - \widehat{T}$ converges to $0.$ 

Now we have to show that $\widehat{L}_T$ satisfies condition $2$ in Lemma \ref{lemma 4.4}, which we are going to prove next, \\

    \begin{claim}
    There exists a $C>0$ such that for any $\delta >0$, there exists $r_\delta>0$ depending on $\epsilon$ such that if $0<r\leq r_\delta,$ then for all $x \in$ spt $r^{-1}\widehat{L}_T \cap B_2(0),$
    $$\frac{1}{\mu_{r^{-1}L_T}(B_\delta(x))}\int_{B_{\delta}(x)} \lambda_r d\mu_{r^{-1}L_T} \geq \frac{1}{2}.$$
    \end{claim}

    \begin{claimproof}
    By Radon-Nikodym theorem, we can find $\lambda_r$ such that 
    $$d\mu_{r^{-1}\widehat{L}_T} = \lambda_r d\mu_{r^{-1}L_T}.$$
    Fix any $\delta >0$. By Neves \cite{neves2007singularities}, there exists a $C$ such that if $x\in spt \;r^{-1}L_T,$ then we have 
    $$\mu_{r^{-1}L_T}\big(B_\delta(x)\big) \geq C_1\delta^2.$$ 
    Also, 
    $$\mu_{r^{-1}L_T} \big(B_1(0)\big) \leq C_2.$$
    As 
    $$\frac{1}{\mu_{r^{-1}L_T} \big(B_1(0)\big)}\int_{B_1(0)} \lambda_r d\mu_{r^{-1}L_T} \to 1,$$ 
    when $r \to 0$, so we can choose $r_0> 0$ such that if $0 <r \leq r_0,$ then
    $$\frac{1}{\mu_{r^{-1}L_T} \big(B_1(0)\big)}\int_{B_1(0)} \lambda_r d\mu_{r^{-1}L_T} \leq 1 - \frac{C_1\delta^2}{4C_2}.$$
    Let,
    $$\frac{1}{\mu_{r^{-1}L_T}(B_\delta(x))}\int_{B_{\delta}(x)} \lambda_r d\mu_{r^{-1}L_T} \geq \frac{1}{2}$$ 
    for some $r,$ such that $0< r \leq r_0$. Then, \\
    \begin{align*}
        1- \frac{C_1\delta^2}{4C_2} &\leq \frac{1}{\mu_{r^{-1}L_T} \big(B_1(0)\big)}\int_{B_1(0)}\lambda_rd\mu_{r^{-1}L_T} \\
        &= \frac{1}{\mu_{r^{-1}L_T}\big(B_1(0)\big)} \bigg(\int_{B_{\delta(x)}} \lambda_rd\mu_{r^{-1}L_T} + \int_{B_1(0) \setminus B_{\delta(x)}}\lambda_rd\mu_{r^{-1}L_T} \bigg) \\
        &\leq \frac{1}{\mu_{r^{-1}L_T}\big(B_1(0)\big)} \bigg(\frac{1}{2} \mu_{r^{-1}L_T}\big(B_\delta(x) \big)  + \int_{B_1(0) \setminus B_{\delta(x)}}\lambda_rd\mu_{r^{-1}L_T} \bigg) \\
        &\leq 1 - \frac{\mu_{r^{-1}L_T}\big(B_\delta(x)\big)}{2\mu_{r^{-1}L_T}(B_1(0))}  \\
        &\leq 1- \frac{C_1\delta^2}{2C_2}.
    \end{align*}
    which implies,
    $$\frac{1}{2} \leq \frac{1}{4},$$
    which is a contradiction. This proves our claim. 
    \end{claimproof} \\

From the last claim, we can conclude that for any $\delta > 0,$ there exists $r_\delta>0$ depending on $\delta$ such that if $0<r\leq r_\delta,$ then for all if $x \in$ spt $r^{-1}\widehat{L}_T \cap B_2(0)$, then 
\begin{align*}
    \mu_{B_\delta (x)} (r^{-1}\widehat{L}_T) &=  \int_{B_{\delta}(x)} \lambda_r d\mu_{r^{-1}L_T} \\
    &\geq \frac{1}{2} \mu_{r^{-1}L_T}(B_\delta(x)) \\  
    &\geq C \delta^2.
\end{align*}
So now we can imitate the proof of the Proposition \ref{prop 4.5} to prove Theorem \ref{thm 3.2}.  We can argue as in the proof to conclude that we can find a sequence $r_n \to \infty $ such that for all $r_n$,
$$M_{B_{r_n}(0)} \widehat{T} \to M_{B_{r_n}(0)} T.$$
But as for all open sets $A$, $\widehat{T}(A) \leq T(A)$, hence by Radon-Nikodym theorem, 
$$\frac{d\mu_{\widehat{T}}}{d\mu_T}  \leq 1$$
$\mu_T$ almost everywhere, but it can not be strictly less than $1$ on a set of positive measure as otherwise it will contradict $M_{B_{r_n}(0)}(T) = M_{B_{r_n}(0)}(\widehat{T})$ for some $r_n$. So we can conclude that 
$$\mu_T = \mu_{\widehat{T}}.$$
This proves the proposition.  
\end{proof}

Now we will argue why proving Theorem \ref{thm 4.1} is enough to prove our main Theorem \ref{thm 1.1}. Note that Theorem \ref{thm 4.1} says that the tangent cone is unique for the current $\widehat{L}_T$. But in Proposition $\ref{prop 4.5},$ we have established the correspondence between the tangent flow of the original flow $L_t$ and the tangent cone of limiting current $\widehat{L}_T$. So by Proposition $\ref{prop 4.5},$ showing the uniqueness of the tangent cone of limiting current $\widehat{L}_T$ will imply the uniqueness of the tangent flow of the limiting varifold $L_T$ and the uniqueness of the tangent flows of the original flow as well, which is our Theorem \ref{thm 1.1}. Note that this proof is under the assumption that for every sufficiently small $r > 0,$ if there exists a $t_r < T$ close to $T$, depending on $r$ such that for $t_r \leq t < T$, there is only one connected component of $L_t \cap B_{4r}(0)$ that intersects $B_r(0)$. 

Now we will show that the above condition can be dropped to prove Theorem \ref{thm 1.1}. First note that there exists a $r >0$ such that for any $t \geq t_r, t_r$ depending on $r$ we can decompose $L_t$ on $B_r(0)$ as 
$$L_t \cap B_r(0):= \cup_{n \in \{1,...,k\}} L^n_t, \hspace{0.5 cm} t \geq t_{r_0},$$
where for any $r <r_0, L^n_t \cap B_{r}(0)$ is connected in $B_{4r}(0).$ This can be shown as follows,

We have shown that for every $r>0,$ there exists a $t < T$ denoted by $t_r$ such that if $L_{t_r} \cap B_r (0)$ has more than one connected component in $B_{4r}(0)$, then they remain disconnected on $B_r(0)$ for all $t_r \leq t < T$. 

Suppose for some $r_1>0, L_t \cap B_{4r_1}(0)$ which intersects $B_{r_1}(0)$ has more than one connected components for $t \geq t_{r_1}.$ Then for some $k_1 >1,$
$$L_t \cap B_{r_1}(0):= \cup_{n \in \{1,...,k_1\}} L^n_t, \hspace{0.5 cm} t \geq t_{r_1}.$$ 
Also, we assumed that each $n \in \{1,...,k_1\}$ has a contribution to the singularity, i.e if any of the components does not go to $0$ as $t \to 0,$ then we can throw it out.

Now for some $r_2 < r_1,$ if $L^n_t \cap B_{4r_2}(0)$ has more than one connected components for $t \geq t_{r_2},$ then
$$L_t \cap B_{r_2}(0):= \cup_{n \in \{1,...,k_2\}} L^n_t, \hspace{0.5 cm} t \geq t_{r_2}.$$ 
where $k_2 > k_1.$ 

Now this process will eventually stop as the number of connected components of $L_t \cap B_{4r}(0)$ in $B_{r}(0)$ has a uniform bound [\citealp{neves2007singularities}, Lemma $7.2$]. Therefore for some $k \geq 1$ for each $n \in \{1,...,k\}$ each component $ L^n_t$ converges to some $L^n_T$ in the sense of varifold as $t \to T.$ By varifold convergence 
$$L_T \cap B_r(0) = \cup_{n \in \{1,...k\}} L^n_T.$$ 
So by Lemma \ref{lemma 3.7}, the tangent flow of $L_t$ and the union of tangent flows generated by $\{L^n_T\}_{n \in \{1,...,k\}}$ are same. But proving Theorem \ref{thm 4.1} implies that for each $n \in K, L^n_T$ has a unique tangent flow, so their union is also unique. This finishes the proof of Theorem \ref{thm 1.1}. In the rest of the paper, we will prove Theorem \ref{thm 4.1}.

Let for some sequence, $\lambda_i \to \infty, \widehat{L}^i_0$ converges to $\widehat{T}$ as a current. As in the proof of Proposition $\ref{prop 4.5}$, we can conclude that for any ball $B$, 
$$M_B{\widehat{L}^i_0} \to M_B \widehat{T}.$$ 
As a consequence $\theta(\widehat{L}_T,0)$ exists and 
$$\theta(\widehat{L}_T,0) = \frac{M_{B_r(0)}(\widehat{T})}{r^2}$$
for all but countably many $r.$ By Neves \cite{neves2007singularities}, $T$ is an area minimizing cone and $2$-dimensional area minimizing cones are all sum of planes intersecting only at the origin. So the flat norm of current $(\widehat{L}^i_T - \widehat{T}) \;\cap \;B_r(0) \to 0$ for every $r > 0$. Now let us define,
$$f(r) := M_{B_r(0)} (\widehat{L}_T) - r^2 \theta.$$ 
Since $r \mapsto M_{B_r(0)} (\widehat{L}_T)$ is increasing, so $f$ is differentiable almost everywhere. By co-Area formula, we have 
$$f'(r) \geq M_{B_r(0)} (\widehat{L}_T) - 2r \theta.$$
We know from Lemma \ref{lemma 4.4} that there exists an $\epsilon_{\widehat{T}} >0$ such that if the flat norm of $r^{-1}\widehat{L}_T - \widehat{T}$ is smaller than $\epsilon_{\widehat{T}}$, then we can apply the result of Lemma \ref{lemma 4.4}. The set of tangent cones of $\widehat{L}_T$ is compact in flat norm topology. This is a standard fact, for proof, one can look at the compactness theorem by Federer [\citealp{federer2014geometric}, Theorem 4.2.17]. So there exists a $r_0 > 0$, $N \in \mathbb{N}$ and $A = \{\cup_{i=1}^N \widehat{T}_i, \widehat{T}_i$ area minimizing cones\} such that if $0 < r < r_0,$ then there exists an area minimizing cone $\widehat{T} \in A$, such that the flat norm of $r^{-1}\widehat{L}_T - \widehat{T}$ is smaller than $\epsilon$, where $\epsilon$ is defined by 
$$\epsilon: = min \{\epsilon_{\widehat{T}_i}: \widehat{T}_i \in A\}.$$
So from Lemma $\ref{lemma 4.4}$,
\begin{align*}
        M_{B_1(0)} (r^{-1} \widehat{L}_T) - r^2\theta &\leq C_r(1- \epsilon)\Big(M_{B_1(0)} \big(0 \; {\RangeX}\;\partial (r^{-1}L_T \mres B_1(0))\big) \\ 
        &\hspace{0.5 cm} - M_{B_1(0)}(\widehat{T}) \Big) +  (C_r-1) M_{B_1(0)} (\widehat{T}). 
    \end{align*}
Hence,    
\begin{align*}
    M_{B_r(0)} (\widehat{L}_T) - r^2\theta  &\leq C_r(1- \epsilon)\Big(M \big(0\; {\RangeX}\; \partial (\widehat{L}_T \cap B_r(0))\big) -  r^2\theta\Big) + (C_{r} -1)r^2\theta\\
    &= C_r (1- \epsilon)\Big( \frac{r}{2} M \big(\partial (\widehat{L}_T \mres B_r(0))\big) -  r^2\theta\Big)+ (C_r -1)r^2 \theta\\
    &= C_r (1- \epsilon)\frac{r}{2}\Big(M \big(\partial (\widehat{L}_T \mres B_r(0))\big) -  2r\theta\Big)+ Cr^4 \theta.
\end{align*}
So, 
$$f(r) \leq C_r(1 - \epsilon) \frac{r}{2} f'(r) + Cr^4 \theta$$
$$f'(r) - \frac{2}{rC_r (1-\epsilon)}f(r) \geq -\frac{2Cr^3\theta}{(C_r(1 -\epsilon))} $$
$$\bigg(\frac{C_r}{r^2}\bigg)^{\frac{1}{1-\epsilon}} f'(r) - \frac{2}{r^3(1-\epsilon)}\bigg(\frac{C_r}{r^2}\bigg)^{\frac{\epsilon}{1-\epsilon}} f(r) \geq  -\frac{2C \theta r}{1-\epsilon}\bigg(\frac{C_r}{r^2} \bigg)^\frac{\epsilon}{1-\epsilon}. \newline$$ \\
Choose $r$ small, so that $C_r \leq 2$, then 
$$-C_r^{\frac{\epsilon}{1-\epsilon}} \leq -2^\frac{\epsilon}{1-\epsilon}.$$ 
Hence,
$$\frac{d}{dr} \bigg( \bigg(\frac{C_r}{r^2}\bigg)^{\frac{1}{1-\epsilon}}  f(r)\bigg) \geq -\overline{C} r^{1-\frac{2\epsilon}{1-\epsilon}}.$$ \\
Integrating from $t$ to $s,$ we get
$$\bigg(\frac{C_s}{s^2}\bigg)^{\frac{1}{1-\epsilon}} f(s) - \bigg(\frac{C_t}{t^2}\bigg)^{\frac{1}{1-\epsilon}}  f(t) \geq - \overline{C} \bigg( s^{\frac{2 - 4\epsilon}{1-\epsilon}} - t^{\frac{2 - 4\epsilon}{1-\epsilon}} \bigg)$$
$$\frac{f(t)}{t^2} \leq \frac{s^2}{t^2} \bigg(\frac{C_s}{s^2}\bigg)^{\frac{1}{1-\epsilon}} \bigg(\frac{C_t}{t^2}\bigg)^{\frac{1}{\epsilon-1}} \frac{f(s)}{s^2} + \bigg(\frac{C_t}{t^2}\bigg)^{\frac{1}{\epsilon - 1}} \overline{C} \bigg( s^{\frac{2 - 4\epsilon}{1-\epsilon}} - t^{\frac{2 - 4\epsilon}{1-\epsilon}} \bigg).$$ 
\\
Now let us denote the function $x \mapsto \frac{x}{|x|}$ by $F$. Then by Area formula for $0 < t < s < \delta,$ we have,
\begin{align*}
    M (F_\# (\widehat{L}_T \mres (B_s(0) \setminus B_t(0)) &\leq \int_{B_s(0) \setminus B_t(0)} \frac{|x^{\perp}|}{|x|^3} \;d \mu_{\widehat{L}_T} \\
    &\leq \bigg(\int_{B_s(0) \setminus B_t(0)} \frac{|x^{\perp}|^2}{|x|^4} \;d \mu_{\widehat{L}_T} \bigg)^{\frac{1}{2}} \bigg(\int_{B_s(0) \setminus B_t(0)} \frac{1}{|x|^2} \;d \mu_{\widehat{L}_T} \bigg)^{\frac{1}{2}} \\
    &\leq \bigg( \frac{f(s)}{\omega_2s^2} - \frac{f(t)}{\omega_2t^2} + \epsilon(s,t) \bigg)^{\frac{1}{2}} \bigg(\frac{M_{B_s(0)} (\widehat{L}_T)}{t^2}\bigg)^\frac{1}{2} \\
    &\leq \bigg(\frac{1}{s^2} \bigg(\frac{C_s}{s^2}\bigg)^{\frac{1}{\epsilon-1}} \Big(g(\delta) \frac{f(\delta)}{\omega_2 \delta^2} + h(\delta) \Big) + \epsilon(s,t)\bigg)^{\frac{1}{2}} \\ 
    &\hspace{0.5 cm}\bigg(\bigg(\frac{s}{t}\bigg)^2 \bigg(\frac{M_{B_\delta(0)} (\widehat{L}_T)}{\delta^2} + \epsilon(\delta,t)\bigg)^{\frac{1}{2}}.
\end{align*}
where, 
$$g(\delta) = \delta^2 \bigg(\frac{C_\delta}{\delta^2}\bigg)^{\frac{1}{1-\epsilon}}, \hspace{0.4 cm}h(\delta) = \overline{C} \delta^{\frac{2-4\epsilon}{1-\epsilon}}.$$ 
Now, if we choose $\eta = \frac{2\epsilon}{1-\epsilon}$ then for small $s$
$$\frac{1}{s^2} \bigg(\frac{C_s}{s^2}\bigg)^{\frac{1}{\epsilon-1}} \leq Cs^\eta.$$
Also,
$$\frac{M_{B_\delta(0)} (\widehat{L}_T)}{\delta^2} + \epsilon(\delta,t) \leq C.$$
Using the above inequalities and Remark $\ref{rmk 4.2}$, we can conclude that if $0 < \frac{s}{2} <t <s,$ then we have, 
$$M(F_\# (\widehat{L}_T \mres (B_s(0) \setminus B_t(0))))\leq Cs^{\frac{\eta}{2}}.$$ 
Hence by iteration on diadic intervals, for every $n$,
$$M (F_\# (\widehat{L}_T \mres (B_s(0) \setminus B_{\frac{s}{2^n}}(0)))) \leq 2Cs^{\frac{\eta}{2}}.$$ 
So this implies that the currents $\partial(\widehat{L}_T \cap B_r(0))$ is cauchy with respect to the flat norm and hence it converges to a unique current $\mathcal{C}.$ Hence $\widehat{L}_T \cap B_r(0)$ converges to the unique cone $0\; {\RangeX}\;\mathcal{C}.$ This proves Theorem \ref{thm 4.1}.

\bibliographystyle{plain}
\bibliography{references.bib}

\end{document}